\documentclass{article}
\usepackage{amsthm,epsfig,a4}
\usepackage[ansinew]{inputenc}
\usepackage{amsmath,amssymb}
\usepackage{amsfonts}
\usepackage{eufrak}
\def\q{\nolinebreak \hfill $\Box$} 
\def\N{\\[5mm] \indent} 

\def\:{\colon\thinspace}
\theoremstyle{plain}
\newtheorem{Lemma}{Lemma}[section]
\newtheorem{Prop}[Lemma]{Proposition}
\newtheorem{Thm}[Lemma]{Theorem}
\newtheorem{Cor}[Lemma]{Corollary}
\newtheorem{Ex}[Lemma]{Example}
\newtheorem{SLemma}{Lemma}[subsection]
\newtheorem{SProp}[SLemma]{Proposition}
\newtheorem{SThm}[SLemma]{Theorem}
\newtheorem{SCor}[SLemma]{Corollary}
\theoremstyle{definition}
\newtheorem{Def}[Lemma]{Definition}

\newtheorem*{Rem}{Remark}
\newtheorem*{Ac}{Acknowledgement}

\begin{document}
\title{Geography of Non-Formal Symplectic and Contact Manifolds}
\author{Christoph Bock}
\date{}
\maketitle
\begin{abstract}
Let $(m,b)$ be a pair of natural numbers. For $m$ even (resp.\ $m$
odd and $b \ge 2$) we show that if there is an $m$-dimensional
non-formal compact oriented manifold with first Betti number $b_1 =
b$, there is also a symplectic (resp.\ contact) manifold with these
properties.
\end{abstract}
\section{Introduction}
In \cite{FMGeo} Fern\'andez and Mu\~noz proved the following:
\begin{Thm} $\,$ \label{nicht sympl}
Given $m \in \mathbb{N}_+$ and $b \in \mathbb{N}$, there are
compact oriented $m$-dimensional manifolds with $b_1 = b$ which
are non-formal if and only if one of the following conditions
holds:
\begin{itemize}
\item[(i)] $m \ge 3 \mbox{ and } b \ge 2$,

\item[(ii)] $m \ge 5 \mbox{ and } b = 1$,

\item[(iii)] $m \ge 7 \mbox{ and } b = 0$. \q
\end{itemize}
\end{Thm}

A natural question to ask is when there are such examples of
symplectic $m$-dimensional manifolds with $b_1 = b$ which are not
formal. Clearly, $m$ has to be even. The same authors as above
constructed in \cite{FMacht} a simply-connected $8$-dimensional
example. Taking products with $S^2$, we get simply-connected
examples for all even dimensions greater than eight.

Fern\'andez, Gotay and Gray showed in \cite{FGG} that there are
$T^2$ bundles over $T^2$ which are symplectic, non-formal with $b_1
= 2,3$ (and of course $4$-dimensional). Again, by taking products
with $S^2$, one also has examples for all even dimensions greater
than four with the same $b_1$.

We prove:
\begin{Thm} $\,$ \label{Main1}
For all $m \in 2\mathbb{N}$ and $b \in \mathbb{N}$, $m,b \ge 4$,
there are compact $m$-dimensional symplectic manifolds with $b_1 =
b$ which are non-formal.
\end{Thm}

\begin{Thm} $\,$ \label{Main2}
For all $m \in 2\mathbb{N}$, $m \ge 6$, there are compact
$m$-dimensional symplectic manifolds with $b_1 = 1$ which are
non-formal.
\end{Thm}

These theorems and the considerations from above imply:
\begin{Thm} $\,$ \label{Fazit}
Let $(m,b) \in 2\mathbb{N}_+ \times \mathbb{N}$. If there is a
non-formal compact oriented $m$-dimensional manifold with $b_1 = b$,
then there is also a symplectic manifold with these properties.  \q
\end{Thm}

Moreover, we will prove:
\begin{Thm} \label{Main}
For each pair $(m,b)$ with $m$ odd and $b \ge 2$ there exists a non-formal
compact contact $m$-manifold.
\end{Thm}
\section{Formality and s-Formality}
We give a brief review of the notion of formality.

A \emph{differential graded algebra (DGA)} is a graded
$\mathbb{R}$-algebra $A = \bigoplus_{i \in \mathbb{N}}~ A^i$
together with an $\mathbb{R}$-linear map $d \: A \to A$ such that
$d(A^i) \subset A^{i+1}$ and the following conditions are satisfied:
\begin{itemize}
\item[(i)] The $\mathbb{R}$-algebra structure of $A$ is given by an
inclusion $\mathbb{R} \hookrightarrow A^0$.

\item[(ii)] The multiplication is graded commutative, i.e.\ for $a \in
A^i$ and $b \in A^j$ one has $a \cdot b = (-1)^{i \cdot j} b \cdot a
\in A^{i+j}$.

\item[(iii)] The Leibniz rule holds: $\forall_{a \in A^i} \forall_{b
\in A} ~ d(a \cdot b) = d(a) \cdot b + (-1)^i a \cdot d(b)$

\item[(iv)] The map $d$ is a differential, i.e.\ $d^2 = 0$.
\end{itemize}
Further, we define $|a| := i$ for $a \in A^i$.

The \emph{$i$-th cohomology of a DGA} $(A,d)$ is the algebra
$$ H^i(A,d) := \frac{\ker (d \: A^i \to A^{i+1})}{\mathrm{im}
(d \: A^{i-1} \to A^i)}.$$

If $(B,d_B)$ is another DGA, then an $\mathbb{R}$-linear map $f \: A
\to B$ is called \emph{morphism} if $f(A^i) \subset B^i$, $f$ is
multiplicative, and $d_B \circ f = f \circ d_A$. Obviously, any such
$f$ induces a homomorphism $f^* \: H^*(A,d_A) \to H^*(B,d_B)$. A
morphism of differential graded algebras inducing an isomorphism on
cohomology is called \emph{quasi-isomorphism}.
\begin{Def} $\,$
A DGA $(\mathcal{M},d)$ is said to be \emph{minimal} if
\begin{itemize}
\item[(i)] there is a graded vector space $V = \big( \bigoplus_{i \in
\mathbb{N}_+} V^i \big) = \mathrm{Span} \, \{a_k ~ | ~ k \in I\}$
with homogeneous elements $a_k$, which we call the generators,

\item[(ii)] $\mathcal{M} = \bigwedge V$,

\item[(iii)]  the index set $I$ is well ordered, such that $k < l ~
\Rightarrow |a_k| \le |a_l|$ and the expression for $da_k$ contains
only generators $a_l$ with $l < k$.
\end{itemize}
\end{Def}
We shall say that $(\mathcal{M},d)$ is a \emph{minimal model for a
differential graded algebra} $(A,d_A)$ if $(\mathcal{M},d)$ is
minimal and there is a quasi-isomorphism 
$\rho \: (\mathcal{M},d) \to (A,d_A)$, i.e.\ $\rho$ induces an
isomorphism $\rho^* \: H^*(\mathcal{M},d) \to H^*(A,d_A)$  on
cohomology.

The importance of minimal models is reflected by the following
theorem, which is taken from Sullivan's work \cite[Section
5]{Sullivan}.
\begin{Thm}\label{Konstruktion d m M}
$\,$ A differential graded algebra $(A,d_A)$ with $H^0(A,d_A) =
\mathbb{R}$ possesses a minimal model. It is unique up to
isomorphism of differential graded algebras.
\end{Thm}

A \emph{minimal model} $(\mathcal{M}_M,d)$ \emph{of a connected
smooth manifold} $M$ is a minimal model for the de Rahm complex
$(\Omega(M),d)$ of differential forms on $M$. The last theorem
implies that every connected smooth manifold possesses a minimal
model which is unique up to isomorphism of differential graded
algebras.

For the remainder of this section, we deal with the notion of
formality. Endowed with the trivial differential, the cohomology of
a minimal DGA is a DGA itself, and therefore it also possesses a
minimal model. In general, these two minimal models need not to be
isomorphic.

A minimal differential graded algebra $(\mathcal{M},d)$ is called
\emph{formal} if there is a morphism of differential graded algebras
$$\psi \: (\mathcal{M},d) \longrightarrow (H^*(\mathcal{M},d),d_H=
 0)$$ that induces the identity on cohomology.

This means that $(\mathcal{M},d)$ and $(H^*(\mathcal{M},d),d_H=
 0)$ share their minimal model. The following theorem gives an
equivalent characterisation.

\begin{Thm}[{\cite[Theorem 1.3.1]{TO}}] $\,$ \label{form}
A minimal model $(\mathcal{M},d)$ is formal if and only if we can
write $\mathcal{M} = \bigwedge V$ and the space $V$ decomposes as a
direct sum $V = C \oplus N$ with $d(C) = 0$, $d$ is injective on
$N$, and such that every closed element in the ideal $I(N)$
generated by $N$ in $\bigwedge V$ is exact. \q
\end{Thm}

This allows us to give a weaker version of the notion of formality.

\begin{Def} $\,$ \label{s-form}
A minimal model $(\mathcal{M},d)$ is called
$\emph{s\mbox{-formal}}$, $s \in \mathbb{N}$, if we can write
$\mathcal{M} = \bigwedge V$ and for each $i \le s$ the space $V^i$
generated by generators of degree $i$ decomposes as a direct sum
$V^i = C^i \oplus N^i$ with $d(C^i) = 0$, $d$ is injective on $N^i$
and such that every closed element in the ideal $I(\bigoplus_{i \le
s} N^i)$ generated by $\bigoplus_{i \le s} N^i$ in $\bigwedge \big(
\bigoplus_{i \le s} V^i \big)$ is exact in $\bigwedge V$.
\end{Def}

Obviously, formality implies $s$-formality for every $s$.

A connected smooth manifold is called \emph{formal} (resp.\
\emph{$s$-formal}) if its minimal model is formal (resp.\
$s$-formal).

We end this section with some results that allow an easier detection
of formality resp.\ non-formality. The next theorem shows the reason
for defining $s$-formality: in certain cases $s$-formality is
sufficient for a manifold to be formal.

\begin{Thm}[{\cite[Theorem 3.1]{FMDon}}] $\,$ \label{formal = n-1 formal}
Let $M$ be a connected and orientable compact smooth manifold of
dimension $2n$ or $(2n-1)$.

Then $M$ is formal if and only if it is $(n-1)$-formal. \q
\end{Thm}

\begin{Ex}[{\cite[Corollary 3.3]{FMDon}}] $\,$ \label{Bsp s-formal}
\begin{itemize}
\item[(i)] Every connected and simply-connected compact smooth
manifold is $2$-formal.

\item[(ii)] Every connected and simply-connected compact smooth
manifold of dimension seven or eight is formal if and only if it is
$3$-formal. \q
\end{itemize}
\end{Ex}

\begin{Prop}[{\cite[Lemma 2.11]{FMDon}}] $\,$ \label{prod}
Let $M_1,M_2$ be connected smooth manifolds. They are both formal
(resp.\ $s$-formal) if and only if $M_1 \times M_2$ is formal
(resp.\ $s$-formal). \q
\end{Prop}
\section{Massey Products}
An important tool for detecting non-formality is the concept of
Massey products: As we shall see below, the triviality of the Massey
products is necessary for formality.

Let $(A,d)$ be a differential graded algebra.
\begin{itemize}
\item[(i)] Let $a_i \in H^{p_i}(A)$, $p_i > 0,~ 1\le i \le 3$,
satisfying $a_j \cdot a_{j+1} = 0$ for $j = 1,2$. Take elements
$\alpha_i$ of $A$ with $a_i = [\alpha_i]$ and write $\alpha_j \cdot
\alpha_{j+1} = d \xi_{j,j+1}$ for $j = 1,2$. The
\emph{(triple-)Massey product}\label{DMP} $\langle a_1,a_2,a_3
\rangle$ of the classes $a_i$ is defined as
$$ [\alpha_1 \cdot \xi_{2,3} +
(-1)^{p_1 +1} \xi_{1,2} \cdot \alpha_3] \in \frac{H^{p_1 + p_2 + p_3
- 1}(A)}{a_1 \cdot H^{p_2 + p_3 -1}(A) + H^{p_1 + p_2 -1}(A) \cdot
a_3 }.$$

\item[(ii)] Now, let $k \ge 4$ and $a_i \in H^{p_i}(A)$, $p_i > 0,~
1\le i \le k$, such that $\langle a_1,\ldots,a_{k-1} \rangle$ and
$\langle a_2,\ldots,a_k \rangle$ are defined and vanish
simultaneously, i.e.\ there are elements $\xi_{i,j}$ of $A$, $1 \le
i \le j \le k$, $(i,j) \ne (1,k)$, such that
\begin{eqnarray}
a_i = [\xi_{i,i}] & \mbox{and} & d \xi_{i,j} = \sum_{l=i}^{j-1}
\overline{\xi}_{i,l} \cdot \xi_{l+1,j}, \label{hoeheres Massey
Produkt}
\end{eqnarray}
where $\overline{\xi} = (-1)^{|\xi|} \xi.$ The \emph{Massey product}
$\langle a_1,\ldots,a_k \rangle$ of the classes $a_i$ is defined as
the set $\{ [ \sum_{l=1}^{k-1} \overline{\xi}_{1,l} \cdot
\xi_{l+1,k} ] \,|\, \xi_{i,j} \mbox{ satisfies (\ref{hoeheres Massey
Produkt})} \} \subset H^{p_1 + \ldots + p_k - (k-2)}(A)$.

We say that $\langle a_1,\ldots,a_k \rangle$ vanishes if $0 \in
\langle a_1,\ldots,a_k \rangle$.
\end{itemize}

\begin{Rem}
The definition of the triple-Massey product in (i) as an element of
a quotient space is well defined, see e.g.\ \cite[Section 1.6]{TO}.
\end{Rem}

The next two lemmata show the relation between formality (resp.\
$s$-formality) and Massey products.

\begin{Lemma}[{\cite[Theorem 1.6.5]{TO}}] $\,$
For any formal minimal differential graded algebra all Massey
products vanish. \q
\end{Lemma}

\begin{Lemma}[{\cite[Lemma 2.9]{FMDon}}] $\,$ \label{not s-formal}
Let $(A,d)$ be an $s$-formal minimal differential graded algebra.
Suppose that there are cohomology classes $a_i \in H^{p_i}(A)$, $p_i
>0$, $1 \le i \le k$, such that $\langle a_1,\ldots,a_k \rangle$ is
defined. If $p_1 + \ldots + p_{k-1} \le s+k-2$ and $p_2 + \ldots +
p_k \le s+k-2$, then $\langle a_1,\ldots,a_k \rangle$ vanishes. \q
\end{Lemma}

In \cite{FMacht}, Fern\'andez and Mu\~noz introduce a different type
of Massey product, called a-Massey product:
\begin{Def} $\,$ \label{a-Massey}
Let $(A,d)$ be a DGA and let $a,b_1,b_2,b_3 \in H^2(A)$ satisfying
$a \cdot b_i = 0$ for $i = 1,2,3$. Take choices of representatives
$a = [\alpha], b_i = [\beta_i]$ and $\alpha \cdot \beta_i = d\xi_i$
for $i=1,2,3$. Then the \emph{a-Massey product} $\langle
a;b_1,b_2,b_3 \rangle$ is defined as $ [ \xi_1 \cdot \xi_2 \cdot
\beta_3 + \xi_2 \cdot \xi_3 \cdot \beta_1 + \xi_3 \cdot \xi_1 \cdot
\beta_2]$ in
$$\frac{H^8(A)}{\langle b_1,a,b_2 \rangle \cdot H^3(A) + \langle
b_1,a,b_3 \rangle \cdot H^3(A) + \langle b_2,a,b_3 \rangle \cdot
H^3(A)}.$$
\end{Def}

\begin{Lemma}[{\cite[Proposition 3.2]{FMacht}}] $\,$
If a minimal differential graded algebra is formal, then every
a-Massey product vanishes. \q
\end{Lemma}

\begin{Cor}\label{NVDR} $\,$
If the de Rahm complex $(\Omega(M),d)$ of a smooth manifold $M$
possesses a non-vanishing Massey or a-Massey product, then $M$ is
not formal.

If there are cohomology classes $a_i \in H^{p_i}(M)$, where $p_i
>0$ and $1 \le i \le k$, with $p_1 + \ldots + p_{k-1} \le s+k-2$
and  $p_2 + \ldots + p_k \le s+k-2$ such that $\langle
a_1,\ldots,a_k \rangle$ does not vanish, then M is not $s$-formal.
\end{Cor}

\textit{Proof.} This holds since a minimal model $\rho \:
(\mathcal{M}_M,d) \to (\Omega(M),d)$ induces an isomorphism on
cohomology. \q
\section{Donaldson Submanifolds}\label{DUmgfen}
Our examples of non-formal symplectic manifolds will be constructed
in a similar way as in the article \cite{FMDon} of Fern\'andez and
Mu\~noz. The examples will be Donaldson submanifolds of non-formal
symplectic manifolds. Therefore, we quote in this section parts of
\cite{FMDon}.

For the remainder of the section we denote the de Rham cohomology of
a smooth manifold $M$ by $H^*(M)$.

In \cite{D} the following is proven: Let $(M,\omega)$ be a
$2n$-dimensional compact symplectic manifold with $[\omega] \in
H^2(M)$ admitting a lift to an integral cohomology class.
Then there exists $k_0 \in \mathbb{N}_+$ such that for each $k \in
\mathbb{N}_+$ with $k \ge k_0$ there is a symplectic submanifold $j
\: Z \hookrightarrow M$ of dimension $2n-2$ whose Poincar\'e dual
satisfies $\mathrm{PD}[Z] = k [\omega]$.
Moreover, the map $j$ is a homology $(n-1)$-equivalence in the
following sense.

Let $f \: M_1 \to M_2$ be a smooth map between smooth manifolds. $f$
is called \emph{homology $s$-equivalence}, $s \in \mathbb{N}$, if it
induces isomorphisms $f^* \: H^i(M_2) \to H^i(M_1)$ on cohomology
for $i \le s-1$ and a monomorphism for $i=s$.

A symplectic submanifold $j \: Z \hookrightarrow M$ as above is
called \emph{symplectic divisor} or \emph{Donaldson submanifold}.

Concerning minimal models and formality in this context, we quote
the following results. Part (i) resp.\ (ii) in the theorem coincides
with Proposition $5.1$ resp.\ Theorem $5.2$ (i) in \cite{FMDon},
where a proof is given.

\begin{Thm}[\cite{FMDon}] $\,$
\label{FMDon5.1f} Let $f \: M_1 \to M_2$ be a homology
$s$-equivalence between connected smooth manifolds. Denote by
$\rho_i \: (\bigwedge V_i ,d) \to (\Omega(M_i),d)$ the minimal
models  of $M_i$ for $i = 1,2$.
\begin{itemize}
\item[(i)] There exist a morphism $F \: (\bigwedge V_2^{\le s} ,d)
\to (\bigwedge V_1^{\le s} ,d)$ of differential graded algebras such
that $F \: V_2^{< s} \to V_1^{< s}$ is an isomorphism,  $F \: V_2^s
\to V_1^s$ is a monomorphism and $\rho_1^* \circ F^* = f^* \circ
\rho_2^*$.

\item[(ii)] If $M_2$ is $(s-1)$-formal, then $M_1$ is $(s-1)$-formal.
\q
\end{itemize}
\end{Thm}

\begin{Cor}[{\cite[Theorem 5.2(ii)]{FMDon}}] $\,$
Let $M$ be a $2n$-dimensional compact symplectic manifold and $j : Z
\hookrightarrow M$ a Donaldson submanifold.

Then for each $s \le n-2$, we have: If $M$ is $s$-formal, then $Z$
is $s$-formal.

In particular, $Z$ is formal if $M$ is $(n-2)$-formal. \q
\end{Cor}

Next, we want to give a criterion for a Donaldson submanifold not to
be formal.

\begin{Prop} $\,$
Let $M$ be a compact symplectic manifold of dimension $2n$, where $n
\ge 3$.
Using the notation from page \pageref{DMP}, we suppose that there are 
cohomology classes $a_i = [\alpha_i] \in H^1(M), 
1 \le i \le 3$, such that the (triple-)Massey product 
$$\langle a_1, a_2, a_3 \rangle = 
[\alpha_1 \wedge \xi_{2,3} +  \xi_{1,2} \wedge \alpha_3] 
\in \frac{H^{2}(M)}{a_1 \cup H^{1}(M) + H^{1}(M) \cup a_3 }$$ 
is defined and does not vanish.

Then every Donaldson submanifold of $M$ is not $1$-formal.
\end{Prop}

\textit{Proof.} Let $j \: Z \hookrightarrow M$ be a Donaldson
submanifold. Since $n \ge 3$, $j$ is a homology $2$-equivalence.
This implies that the (triple-)Massey product
$$\langle j^*a_1, j^*a_2, j^*a_3 \rangle = 
[j^*\alpha_1 \wedge j^*\xi_{2,3} + j^*\xi_{1,2} \wedge j^*\alpha_3] 
\in \frac{H^{2}(Z)}{j^*a_1 \cup H^{1}(Z) + H^{1}(Z) \cup j^*a_3 }$$ 
is defined and does not vanish. Now, Corollary \ref{NVDR} implies
that $Z$ is not $1$-formal. \q
\N As an immediate consequence of the proposition and its proof we get:
\begin{Cor} \label{Donaldson 1-formal}
Let $Z_1,\ldots,Z_k,M$ be compact symplectic manifolds and assume
that $Z_i \hookrightarrow Z_{i+1}$ and $Z_k \hookrightarrow M$ are
Donaldson submanifolds for $i = 1,\ldots,k-1$. We 
suppose that there are cohomology classes $a_i = [\alpha_i] \in H^1(M), 
1 \le i \le 3$, such that the (triple-)Massey product 
$$\langle a_1, a_2, a_3 \rangle = 
[\alpha_1 \wedge \xi_{2,3} +  \xi_{1,2} \wedge \alpha_3] 
\in \frac{H^{2}(M)}{a_1 \cup H^{1}(M) + H^{1}(M) \cup a_3 }$$ 
is defined and does not vanish.

If $\dim Z_1 \ge 4$, then $Z_1$ is not $1$-formal. \q
\end{Cor}

The next lemma will be needed in the proof of Theorem \ref{Main2}.
The proof is taken word by word from the proof of Formula ($5$) in
\cite{FMDon}. 

\begin{Lemma}\label{Don inj} $\,$
Let $(M,\omega)$ be a $2n$-dimensional compact symplectic manifold
and $j \: Z \hookrightarrow M$ a Donaldson submanifold.

Then for each $p = 2(n-1) - i$, $0 \le i \le (n-2)$, there is a
monomorphism
$$ \frac{H^p(M)}{\ker ( [\omega] \cup \: H^p(M) \to H^{p+2}(M) )}
\longrightarrow H^p(Z).$$
\end{Lemma}

\textit{Proof.} The claim can be seen via Poincar\'e duality. Let $0
\le i \le (n-2)$, $p = 2(n-1) - i$ and $\alpha \in \Omega^p(M)$ be
closed. Then we have
$$j^* [\alpha] = 0 \Longleftrightarrow  \forall_{b
\in H^i(Z)}~ j^*[\alpha] \cup b = 0.$$ Since $i \le (n-2)$, we know
that there is an isomorphism $ j^* \: H^i(M) \stackrel{\simeq}{\to}
H^i(Z)$, thus we can assume that for each $b \in H^i(Z)$ there is a
closed $i$-form $\beta$ on $M$ with $[\beta |_Z] = j^*[\beta] = b$
and get
$$ j^*[\alpha] \cup j^*[\beta] = \int_Z j^*\alpha \wedge j^*\beta =
\int_M \alpha \wedge \beta \wedge k \omega, $$ since $[Z] =
\mathrm{PD}[k \omega]$ with $k \in \mathbb{N}_+$. Therefore, we have
$$j^*[\alpha] = 0 \Longleftrightarrow  \forall_{[\beta] \in H^i(M)}~
[\alpha \wedge \omega] \cup [\beta] = 0  \Longleftrightarrow [\alpha
\wedge \omega] =0,$$ from where the lemma follows. \q
\section{Known Examples}
\subsection{The manifolds M(p,q)} \label{M(p,q)}
The following examples are taken from \cite{CFG}.

Let $R$ be a ring with $1$. For $p \in  \mathbb{N}_+$ let $H(1,p;R)$
be the set
$$ \{ \left( \begin{array}{ccc} I_p & x & z \\
0 & 1 & y \\ 0 & 0 & 1 \end{array} \right) \,|\, x,z \in R^p , 
\, y \in R \}.$$ We write $H(1,p)$ for $H(1,p;\mathbb{R})$.
Clearly, this is a nilpotent Lie group and the $2p+1$ differential
$1$-forms
$$ \alpha_i := dx_i,~ \beta := dy,~ \gamma_i :=dz_i - x_i \, dy,~~
1 \le i \le p $$ form a basis of the left-invariant $1$-forms.
Obviously, we have $d\alpha_i = d \beta = 0$ and $d \gamma_i = -
\alpha_i \wedge \beta$.

Further, let $q \in \mathbb{N}_+$. We set $G(p,q) := H(1,p) \times
H(1,q)$. Again, this is a Lie group and analogous as above, we
denote the $2p+2q+2$ forms which form a basis of the left-invariant
$1$-forms by
$$ \alpha_1, \ldots, \alpha_p, \beta, \gamma_1, \ldots, \gamma_p,
\tilde{\alpha}_1, \ldots, \tilde{\alpha}_q, \tilde{\beta},
\tilde{\gamma}_1, \ldots, \tilde{\gamma}_q.$$ An easy computation
shows that the $2$-form $$\omega := \sum_{i=1}^p \alpha_i \wedge
\gamma_i + \sum_{i=1}^q \tilde{\alpha}_i \wedge \tilde{\gamma}_i +
\beta \wedge \tilde{\beta}$$ is a left-invariant symplectic form.
Therefore $M(p,q) := G(p,q) / \Gamma(p,q)$, where $\Gamma(p,q) :=
H(1,p;\mathbb{Z}) \times H(1,q;\mathbb{Z})$, is a compact symplectic
nilmanifold of dimension $2p+2q+2$.

By \cite[Theorem 2.1.3]{TO}, the minimal model $\rho \: \big(
\mathcal{M}_{M(p,q)} , d \big) \to \big( \Omega(M(p,q)),d \big)$ is
given by
\begin{eqnarray}
& \mathcal{M}_{M(p,q)} = \bigwedge (a_1, \ldots, a_p, b, c_1,
\ldots, c_p, \tilde{a}_1, \ldots, \tilde{a}_q, \tilde{b},
\tilde{c}_1, \ldots, \tilde{c}_q ), & \nonumber \\
& |a_i| = |b| = |c_i| = |\tilde{a}_i| = |\tilde{b}| = |\tilde{c}_i|
= 1, & 
\nonumber \\ & da_i = db = d\tilde{a}_i =
d\tilde{b} = 0,~ dc_i = -a_i b,~ d\tilde{c}_i = -\tilde{a}_i
\tilde{b}, & \nonumber \\ & \rho(a_i) = \alpha_i,~ \rho(b) = \beta,~
\rho(c_i) = \gamma_i,~ \rho(\tilde{a}_i) = \tilde{\alpha}_i,~
\rho(\tilde{b}) = \tilde{\beta},~ \rho(\tilde{c}_i) =
\tilde{\gamma}_i. & \nonumber
\end{eqnarray}
Therefore, we see $b_1(M(p,q)) = p+q+2$.

\begin{Prop}[\cite{CFG}]\label{Mpqnichtformal} $\,$
$M(p,q)$ is not formal.
\end{Prop}

\textit{Proof.} $\langle [\beta], [\alpha_i], [\alpha_i] \rangle$ is
a non-vanishing Massey product. \q
\N Using Theorem \cite[Theorem 2.1.3]{TO} again, one computes the
first cohomology groups of $M(p,q)$ as

\parbox{11cm}{\begin{eqnarray*}
H^0(M(p,q)) &=& \langle 1 \rangle, \\ H^1(M(p,q)) &=& \langle
[\alpha_i], [\beta], [\tilde{\alpha}_k] , [\tilde{\beta}]
\,|\, 1 \le i \le p, 1 \le k \le q \rangle, \\
H^2(M(p,q)) &=& \langle [\alpha_i \wedge \gamma_j], [\alpha_i \wedge
\tilde{\alpha}_k], [\alpha_i \wedge \tilde{\beta}], [\beta \wedge
\gamma_j], [\beta \wedge \tilde{\alpha}_l], [\beta \wedge
\tilde{\beta}], \\ && ~ [\tilde{\alpha}_k \wedge \tilde{\gamma}_l],
[\tilde{\beta} \wedge \tilde{\gamma}_l] \,|\, 1 \le i,j \le p, 1 \le
k,l \le q \rangle.
\end{eqnarray*}}
%
\subsection{The manifold M$_{8,0}$} \label{M8}
Fern\'andez and Mu\~noz constructed in \cite{FMacht} an
$8$-dimensional non-$3$-formal compact symplectic manifold
$(M_{8,0},\omega)$ with

\parbox{11cm}{\begin{eqnarray*}
& b_0(M_{8,0}) = b_8(M_{8,0}) = 1,~  b_1(M_{8,0}) = b_7(M_{8,0}) = 0
,  &  \\ & b_2(M_{8,0}) = b_6(M_{8,0}) = 256,~ b_3(M_{8,0}) =
b_5(M_{8,0}) = 0,~ b_4(M_{8,0}) = 269 &
\end{eqnarray*}}
\parbox{11mm}{\begin{eqnarray} \label{Betti} \end{eqnarray}}

\noindent as desingularisation of an orbifold. The latter is a
$\mathbb{Z}_3$-quotient of a nilmanifold. The non-formality is
proved by regarding the a-Massey product $\langle [\vartheta];
[\tau_1], [\tau_2], [\tau_3] \rangle$ for certain closed $2$-forms
$\vartheta, \tau_i$ on $M_{8,0}$: One has $\langle [\vartheta];
[\tau_1], [\tau_2], [\tau_3] \rangle = \lambda \,[ \omega^4 ]$ for
$\lambda \ne 0$. Clearly, $\lambda \, \omega^4$ is not exact, and
since $b_3(M_{8,0})=0$, it follows from Definition \ref{a-Massey}
that this a-Massey product does not vanish.
\section{Proofs}
\subsection{Proof of Theorem \ref{Main1}} Because products with
finitely many copies of $S^2$ give the higher-dimensional examples,
it is enough to prove that for every $b \ge 4$ there is a non-formal
compact symplectic $4$-manifold $M$ with $b_1(M) = b$.

Let $b \ge 4$ and choose $p,q \in \mathbb{N}_+$ such that $p+q+2 =
b$. Then $M(p,q)$ has dimension $2p+2q+2 \ge 6$ and is a non-formal
compact symplectic nilmanifold with $b_1(M(p,q)) = b$ 
which satisfies the assumption of Corollary \ref{Donaldson
1-formal}. Therefore, we get the required non-formal $4$-manifold 
$Z$ with $b_1(Z)= b_1(M(p,q)) = b$. \q
\subsection{Proof of Theorem \ref{Main2}} Since direct products with
finitely many copies of $S^2$ gives the higher-dimensional ones, it
is enough to find a six-dimensional example. This is constructed in
\cite{ipse}.

But using the ideas from above, one can construct an eight
dimensional example as follows:

Gompf has shown in \cite{Go} that there is a compact symplectic
$4$-manifold $M_{4,1}$ with $b_1(M_{4,1}) = 1$. By Proposition
\ref{prod}, $M_{12,1} := M_{8,0} \times M_{4,1}$ is a compact
symplectic $12$-manifold which is not $3$-formal. Clearly, we have
$b_1(M_{12,1}) = 1$. Denote the projections by $\pi \: M_{12,1} \to
M_{8,0}$, $\mathfrak{p} \: M_{12,1} \to M_{4,1}$ and the symplectic
forms of $M_{8,0}, M_{4,1}$ and $M_{12,1}$ by $\omega, \sigma$ and
$\Omega = \pi^* \omega + \mathfrak{p}^*\sigma$. Let $\vartheta,
\tau_i$ be the $2$-forms of Section \ref{M8}. We mentioned $\langle
[\vartheta]; [\tau_1], [\tau_2], [\tau_3] \rangle = \lambda \,[
\omega^4] \ne 0$.

Let $j \: Z_{10,1} \hookrightarrow M_{12,1}$ be a Donaldson
submanifold. The $10$-form
\begin{equation*}
\Omega \wedge \lambda\, \pi^* \omega^4 = (\pi^* \omega +
\mathfrak{p}^* \sigma) \wedge \lambda\, \pi^* \omega^4 = \lambda\,
\mathfrak{p}^* \sigma \wedge \pi^* \omega^4
\end{equation*}
on $M_{12,1}$ does not represent the zero class in
$$H^{10}(M_{12,1}) \stackrel{(\ref{Betti})}{=} ( \langle [\sigma^2] \rangle \otimes
H^6(M_{8,0})) \oplus ( H^2(M_{4,1}) \otimes H^8(M_{8,0}) ).$$
Therefore, we get from Lemma \ref{Don inj}: $\lambda\, j^* \pi^*
[\omega^4] \in H^8(Z_{10,1}) \backslash \{0\}$. From (\ref{Betti})
we know $H^5(M_{8,0}) = 0$. Hence $\langle [j^* \pi^* \tau_k], [j^*
\pi^* \vartheta], [j^* \pi^* \tau_l] \rangle = 0$ for $1\le k,l
\le3$. So in the following a-Massey product there is no
indeterminacy:
$$\langle [j^* \pi^* \vartheta]; [j^* \pi^* \tau_1], [j^* \pi^*
\tau_2], [j^* \pi^* \tau_3] \rangle = \lambda\, j^* \pi^* [\omega^4]
\ne 0.$$ It follows that $Z_{10,1}$ is not formal. The fact that
$\dim Z_{10,1} = 10$ and $b_1(Z_{10,1}) = 1$ is clear by the remarks
in Section \ref{DUmgfen}.

Now, let $\tilde{j} \: Z_{8,1} \hookrightarrow Z_{10,1}$ be a
Donaldson submanifold. Then the $10$-form $j^* \Omega \wedge
\lambda\, j^* \pi^* \omega^4$ on $Z_{10,1}$ does not represent the
zero class in $H^{10}(Z_{10,1})$, for we have $$\Omega^2 \wedge
\pi^* \omega^4 = (\mathfrak{p}^* \sigma + \pi^* \omega) \wedge
(\mathfrak{p}^* \sigma \wedge \pi^* \omega^4) = 2\, \mathfrak{p}^*
\sigma^2 \wedge \pi^* \omega^4 \ne 0,$$ and by Lemma \ref{Don inj}
we get $[j^* ( \Omega \wedge \pi^* \omega^4) ] \ne 0$.

Again we use Lemma \ref{Don inj} to see $\lambda \, \tilde{j}^* j^*
\pi^* [\omega^4] \in H^8(Z_{8,1}) \backslash \{0\}$ and can prove
similarly as for $Z_{10,1}$ that $Z_{8,1}$ is not formal. Moreover,
$Z_{8,1}$ is $8$-dimensional and has first Betti number equal to
one. \q
\begin{Rem}
A Donaldson submanifold $Z_{6,1}$ of the manifold $Z_{8,1}$ that we
have constructed in the last proof is formal: From the $2$-formality
of $M_{12,1} = M_{8,0} \times M_{4,1}$ it follows that $Z_{6,1}$ is
$2$-formal and therefore formal by Theorem \ref{formal = n-1
formal}.
\end{Rem}
\subsection{Proof of Theorem \ref{Main}}
Our starting point is a non-formal symplectic manifold. Boothby and
Wang proved that there is a contact manifold which fibres over it
with fibre a circle.

\begin{SThm}[{\cite[Theorem 3]{BW}}] \label{BW}
If $(M,\omega)$ is a compact symplectic manifold whose symplectic
form determines an integral cohomology class of $M$, then the
principal circle bundle $\pi \: E \to M$ with first Chern class
$c_1(\pi) = [\omega] \in H^2(M,\mathbb{Z})$ admits a connection
$\mathrm{1}$-form $\alpha$ such that $\pi^* \omega = d \alpha$  and
$\alpha$ is a contact form on $E$. \q
\end{SThm}

Let $E,M$ be as in the last theorem. Since $E$ is an $S^1$-bundle
over $M$, one can apply the Gysin sequence to obtain $b_1(E) =
b_1(M)$. We can even find a contact manifold which has the same
fundamental group as $M$:

\begin{SCor} \label{piE=piM}
Let $(M,\omega)$ be a compact symplectic manifold of dimension $2n$
whose symplectic form determines an integral cohomology class.

Then there is a compact contact manifold $(E,\ker \alpha)$ and a
principal circle bundle $\pi \: E \to M \#
\overline{\mathbb{C}P}{^n}$ with first Chern class $c_1(\pi) =
[\omega]$ such that the fundamental groups satisfy $\pi_1(E) =
\pi_1(M \# \overline{\mathbb{C}P}{^n}) = \pi_1(M)$.
\end{SCor}

\textit{Proof.} We use the same argumentation as in the proof of
\cite[Theorem 4.4]{Go}. After blowing up a point in $M$, we can
obtain a manifold $M^{\prime} := M \# \overline{\mathbb{C}P}{^n}$
with a symplectic form $\omega^{\prime}$ such that
$[\omega^{\prime}] = [\omega] + \varepsilon e \in H^2(M^{\prime}) =
H^2(M) \oplus H^2(\overline{\mathbb{C}P}{^n})$, where $\varepsilon
\in \frac{1}{\mathbb{N}_+}$ and $e$ is a generator of
$H^*(\overline{\mathbb{C}P}{^n})$. Without loss of generality, we
can assume that $\omega^{\prime}$ determines an integral cohomology
class and there is an embedded sphere $S \subset M^{\prime} = M \#
\overline{\mathbb{C}P}{^n}$ such that $\int_S \omega^{\prime} = 1$.
(Since $\int_S \omega^{\prime}$ depends on the size of the ball
removed from $M$ in the blow-up, we may have to enlarge $\omega$ by
an integer scale first.) Let $\pi \: E \to M^{\prime}$ with
$c_1(\pi) = [\omega^{\prime}]$ as in Theorem \ref{BW}. Then the
restriction of the fibration $\pi$ to $S$ is the Hopf fibration,
i.e.\ $\pi^{-1} (S) = S^3$ and the middle map in the following part
of the homotopy sequence is an isomorphism:
$$ \{0\} = \pi_2(\pi^{-1} (S)) \longrightarrow \pi_2(S) \longrightarrow
\pi_1(S^1) \longrightarrow \pi_1(\pi^{-1} (S)) = \{1\}. $$ $\pi_2(S)
\to \pi_1(S^1)$ is an isomorphism. From $S \subset M^{\prime}$ we
get in the following part of the homotopy sequence of the fibration
$\pi$ that the first map is surjective:
$$ \pi_2(M^{\prime}) \longrightarrow \pi_1(S^1) \longrightarrow \pi_1(E)
\longrightarrow \pi_1(M^{\prime}) \longrightarrow \pi_0(S^1) =
\{1\}.$$ This yields an isomorphism $\pi_* \: \pi_1(E) \to
\pi_1(M^{\prime}) = \pi_1(M)$. \q
\N Under certain conditions we can show that the contact manifold
that we have just constructed is not formal.

\begin{SProp} \label{E}
Let $(M,\omega)$ be a compact symplectic manifold of dimension $2n
\ge 4$ whose symplectic form determines an integral cohomology
class. Further, suppose that there are cohomology classes $a_i \in
H^{1}(M)$, $1\le i \le 3$, such that $\langle a_1,a_2,a_3 \rangle$
is a non-vanishing Massey product in $M$.

Then the manifold $E$ of Corollary \ref{piE=piM} is not formal.
\end{SProp}

\textit{Proof.} Let $\pi \: E \to M^{\prime} := M \#
\overline{\mathbb{C}P}{^n}$ be as in Corollary \ref{piE=piM} and the
non-vanishing Massey product $\langle a_1,a_2,a_3 \rangle$ be
defined by a $\mathrm{2}$-form $\alpha_1 \cdot \xi_{2,3} + \xi_{1,2}
\cdot \alpha_3$. (Here we use the notation from page \pageref{DMP}.)
We show:
\begin{eqnarray}
& \pi^* \: H^1(M^{\prime}) \to H^1(E) \mbox{ is an isomorphism.} \label{pi iso} &\\
& 
H^2(M) \cap \ker \big(\pi^* \: H^2(M^{\prime}) \to H^2(E)\big) =
\{0\} \label{Kern pi} &
\end{eqnarray}
Then $\pi^* \alpha_1 \cdot \pi^* \xi_{2,3} + \pi^* \xi_{1,2} \cdot
\pi^*  \alpha_3$ defines the non-vanishing Massey product
\begin{eqnarray*}
\langle \pi^* a_1, \pi^* a_2, \pi^* a_3 \rangle & \in & \frac{\pi^*
\big( H^2(M^{\prime}) \big)}{\pi^* a_1 \cdot \pi^*
\big(H^1(M^{\prime})\big)
+ \pi^* \big(H^1(M^{\prime})\big) \cdot \pi^* a_3} \\
& \subset & \frac{H^2(E)}{\pi^* a_1 \cdot H^1(E) + H^1(E) \cdot
\pi^* a_3},
\end{eqnarray*}
so $E$ is not formal.

(Assume $\langle \pi^* a_1, \pi^* a_2, \pi^* a_3 \rangle$ vanishes.
Then for $j=1,2$ there exists a class $[\Xi_{j,j+1}] \in H^1(E)$
such that $0 = d \Xi_{j,j+1} = \pi^* \alpha_j \cdot \pi^*
\alpha_{j+1}$. Property (\ref{pi iso}) implies the existence of
$[\xi_{j,j+1}] \in H^1(M^{\prime})$ with $0 = d \pi^* \xi_{j,j+1} =
\pi^* \alpha_j \cdot \pi^* \alpha_{j+1}$ for $j = 1,2$, i.e.\
$\alpha_j \cdot \alpha_{j+1}$ is exact by (\ref{Kern pi}) and
$\langle [\alpha_1] , [\alpha_2], [\alpha_3] \rangle$ vanishes,
which is a contradiction.)

It remains to show (\ref{pi iso}) and (\ref{Kern pi}): Consider the
Gysin sequence of $\pi$.
\begin{equation} \label{Gysin}
\{0\} \longrightarrow H^1(M^{\prime})
\stackrel{\pi^*}{\longrightarrow} H^1(E) \longrightarrow
H^0(M^{\prime}) \stackrel{[\omega^{\prime}] \cup}{\longrightarrow}
H^2(M^{\prime}) \stackrel{\pi^*}{\longrightarrow} H^2(E)
\longrightarrow \ldots
\end{equation}
Since $\cup [\omega] \: H^0(M^{\prime}) \to H^2(M^{\prime})$ is
injective, it follows that $\pi^* \: H^1(M^{\prime}) \to H^1(E)$ is
an isomorphism.

Further, we get $\ker \big(\pi^* \: H^2(M^{\prime}) \to H^2(E)\big)
\stackrel{(\ref{Gysin})}{=} \mathbb{R} [\omega^{\prime}]$. Denote
$$\mathfrak{pr}_2 \: H^2(M^{\prime}) = H^2(M) \oplus
H^2(\overline{\mathbb{C}P}{^n}) \longrightarrow
H^2(\overline{\mathbb{C}P}{^n})$$ the projection onto the second
factor. Since $\omega^{\prime}$ is the symplectic form of the
blow-up of $M$, we have $\mathfrak{pr}_2 ([\omega^{\prime}]) \ne 0$.
But $\mathfrak{pr}_2|_{H^2(M)} = 0$, so (\ref{Kern pi}) follows. \q
\N Using the preparations that we have done, we are able to
construct explicit non-formal contact manifolds.

\begin{SThm} \label{MainC1}
For each $n \in \mathbb{N}$ with $n \ge 2$ and $b \in \{2,3\}$ there
exists a compact contact $(2n+1)$-manifold which is not formal.
\end{SThm}

\textit{Proof.} In \cite{FGG} the following manifolds are studied.
Let $M_b$, $b \in \{2,3\}$ be the four-dimensional nilmanifold with
basis of left-invariant $1$-forms $\{\alpha, \beta, \gamma, \eta\}$
such that
\begin{eqnarray*}
& d \alpha = d \beta = 0, & \\
& d \gamma = \alpha \wedge \beta, & \\
& d \eta = \left\{
\begin{array}{c@{\quad \mbox{if }}l}
\alpha \wedge \gamma, &  b = 2, \\ 0, & b = 3.
\end{array}
\right.
\end{eqnarray*}
Then, $b_1(M_b) = b$, the $\mathrm{2}$-form $\alpha \wedge \eta +
\beta \wedge \gamma$ is a symplectic form for $M_b$, and $\langle
[\beta], [\beta], [\alpha] \rangle = -[\beta \wedge \gamma]$ is a
non-vanishing Massey product. The case $n=2$ now follows from
Proposition \ref{E}. For $n>2$ consider the manifolds $M_b \times
(S^2)^{n-2}$ instead of $M_b$. \q

\begin{SThm} \label{MainC2}
For each $b \in \mathbb{N}$ with $b \ge 2$, there are non-formal
compact contact manifolds of dimension $3$ and $5$ with first Betti
number $b_1 = b$.
\end{SThm}

\textit{Proof.} By Theorem \ref{nicht sympl}, we know that there is
a compact oriented $3$-manifold $M$ with $b_1 = b \ge 2$ which is
non-formal. By theorems of Martinet \cite{M} and Geiges
\cite[Proposition 2]{Ge} $M$ and $M \times S^2$ admit contact
structures. Further, it follows from Proposition \ref{prod} that $M
\times S^2$ is not formal. \q
\N Now, Theorem \ref{Main} follows from Theorems \ref{MainC1},
\ref{MainC2}, Proposition \ref{prod} and the following result of
Bourgeois:

\begin{SThm}[\cite{Bou}] \label{M x Torus}
Let $M$ be a compact contact manifold of dimension greater than or
equal to three.

Then $M \times T^2$ admits a contact structure. \q
\end{SThm}

Note that the case $(2n+1,b)$ in 
Theorem \ref{Main} is realized if $(2n-1,b-2)$ is realized. Inductively,
one gets to either the case $b \in \{2,3\}$, $2n+1\geq 3$ or the case $b\geq 4$,
$2n+1=3$, both covered previously. \q

\begin{Ac}
The results presented in this paper were obtained in my dissertation
under the supervision of Prof.\ H.\ Geiges. I wish to express my
sincerest gratitude for his support during the last four years.
\end{Ac}

\noindent \author{Christoph Bock \\ Universit\"at zu K\"oln \\
Mathematisches Institut \\ Weyertal 86--90 \\
D-50931 K\"oln \\ e-mail: bock@math.uni-koeln.de}

\end{document}